\magnification=1200
\overfullrule=0pt
\centerline {\bf Integral functionals on $L^p$-spaces: infima over
sub-level sets}\par
\bigskip
\bigskip
\centerline {BIAGIO RICCERI}\par
\bigskip
\bigskip
\centerline {\it Dedicated to Professor Alfonso Villani, with esteem and friendship, on
his sixtieth birthday}\par
\bigskip
\bigskip
\noindent
{\bf Abstract:}  In this paper, we establish the following result:\par 
Let $(T,{\cal F},\mu)$ be a $\sigma$-finite measure space, let
$Y$ be a reflexive real Banach space, and let $\varphi, \psi:Y\to {\bf R}$
be two sequentially weakly lower semicontinuous functionals such
that
$$\inf_{y\in Y}{{\min\{\varphi(y),\psi(y)\}}\over
{1+\|y\|^p}}>-\infty$$
for some $p>0$.
 Moreover,
assume that $\varphi$ has no global minima, while $\varphi+\lambda\psi$
is coercive and has a unique global minimum for each $\lambda>0$.\par
Then, for each
$\gamma\in L^{\infty}(T)\cap L^1(T)\setminus \{0\}$, with $\gamma\geq 0$,
and for each $r>\inf_{Y}\psi$, if we put
$$V_{\gamma,r}=
\left \{u\in L^p(T,Y) : \int_T\gamma(t)\psi(u(t))d\mu\leq
r\int_T\gamma(t)d\mu\right
\}\ ,$$
we have
$$\inf_{u\in V_{\gamma,r}}
\int_T\gamma(t)\varphi(u(t))d\mu=
\inf_{\psi^{-1}(r)}\varphi\int_T\gamma(t)d\mu\ .$$
\bigskip
\noindent
{\bf Key words:} Integral functional; $L^p$-space; coercivity; uniqueness;
global minimum.\par
\bigskip
\noindent
{\bf 2010 Mathematics Subject Classification:} 49J99; 49K99. 
\bigskip
\bigskip
Here and in the sequel, $(T,{\cal F},\mu)$ ($\mu(T)>0$) is a $\sigma$-finite
measure space, $Y$ is a reflexive real Banach space and
$\varphi, \psi:Y\to {\bf R}$ are two sequentially weakly lower semicontinuous
functionals such that
$$\inf_{y\in Y}{{\min\{\varphi(y),\psi(y)\}}\over
{1+\|y\|^p}}>-\infty \eqno{(1)}$$
for some $p>0$. \par
\smallskip
For each $\lambda\in [0,\infty]$, we denote by $M_{\lambda}$ the set of all
global minima of $\varphi+\lambda\psi$ or the empty set according to whether
$\lambda<+\infty$ or $\lambda=+\infty$. We adopt the conventions $\inf\emptyset=+\infty$
and $\sup\emptyset=-\infty$.\par
\smallskip
Moreover, $a, b$ are two fixed numbers in $[0,+\infty]$, with $a<b$, and
$\alpha$, $\beta$ are the
numbers so defined:
$$\alpha=\max\left\{\inf_Y\psi,\sup_{M_b}\psi\right\}\ ,$$
$$\beta=\min\left\{\sup_Y\psi,\inf_{M_a}\psi\right\}\ .$$
\smallskip
As usual, $L^p(T,Y)$ denotes
the space of all $\mu$-strongly measurable functions $u:T\to Y$ such that
$$\int_T\|u(t)\|^pd\mu<+\infty\ .$$
A functional $P:Y\to {\bf R}$ is said to be coercive provided
$$\lim_{\|y\|\to +\infty}P(y)=+\infty\ .$$
The aim of this paper is to establish the following result:\par
\medskip
THEOREM 1. - {\it  Assume that the functional
 $\varphi+\lambda\psi$
is coercive and has a unique global minimum for each $\lambda\in ]a,b[$.
 Assume also that
$$\alpha<\beta\ .$$\par
Then, for each
$\gamma\in L^{\infty}(T)\cap L^1(T)\setminus \{0\}$, with $\gamma\geq 0$,
and for each $r\in ]\alpha,\beta[$, if we put
$$V_{\gamma,r}=
\left \{u\in L^p(T,Y) : \int_T\gamma(t)\psi(u(t))d\mu\leq
r\int_T\gamma(t)d\mu\right
\}\ ,$$
we have
$$\inf_{u\in V_{\gamma,r}}
\int_T\gamma(t)\varphi(u(t))d\mu=
\inf_{\psi^{-1}(r)}\varphi\int_T\gamma(t)d\mu\ .\eqno{(2)}$$}
\medskip
The proof of Theorem 1 is entirely based on the following result that we have
 established in [5]:\par
\medskip
THEOREM A. - {\it Under the assumptions of Theorem 1,
for each $r\in ]\alpha,\beta[$, there exists
$\lambda_r\in ]a,b[$
such that the unique global minimum of $\varphi+\lambda_r \psi$ lies in
$\psi^{-1}(r)$.} \par
\medskip
{\it Proof of Theorem 1}. First, we also assume that
$$\varphi(0)=\psi(0)=0\ .$$
Actually, once we prove the theorem under this additional assumption, 
the general version is obtained applying the particular version to the
functions $\varphi-\varphi(0)$ and $\psi-\psi(0)$.
Next, observe that the functionals
$\varphi$ and $\psi$ are Borel (in the weak topology, and so in the
strong one too). This implies that, for each $u\in L^p(T,Y)$,
the functions $\varphi\circ u$ and $\psi\circ u$ are $\mu$-measurable. On the
other hand,
in view of $(1)$, for some $c>0$, we have
$$-c\gamma(t)(1+\|u(t)\|^p)\leq \gamma(t)\min\{\varphi(u(t)),\psi(u(t))\}$$
for all $t\in T$. Since $\gamma\in L^{\infty}(T)\cap L^1(T)$,
the function $t\to -\gamma(t)(1+\|u(t)\|^p)$ lies in $L^1(T)$, and so
the integrals $\int_T\gamma(t)\varphi(u(t))d\mu$ and
$\int_T\gamma(t)\psi(u(t))d\mu$ exist and belong to $]-\infty,+\infty]$.
For each $\lambda\in ]a,b[$,
denote by $\hat y_{\lambda}$
the unique global minimum in $Y$ of the functional
$\varphi+\lambda\psi$.
 By Theorem A, 
there exists $\lambda_r\in ]a,b[$ such that
$$\psi(\hat y_{\lambda_r})=r\ .$$
So, we have
$$\varphi(\hat y_{\lambda_r})+\lambda_rr\leq \varphi(y)+\lambda_r\psi(y)$$
for all $y\in Y$. From this, it clearly follows that
$$\varphi(\hat y_{\lambda_r})=\inf_{\psi^{-1}(r)}\varphi\ . \eqno{(3)}$$
 Likewise,
for each $u\in L^p(T,Y)$, it follows that
$$(\varphi(\hat y_{\lambda_r})+\lambda_rr)\int_T\gamma(t)d\mu\leq
\int_T(\gamma(t)(\varphi(u(t))+\lambda_r\psi(u(t)))d\mu\ .$$
Therefore, for each $u\in V_{\gamma,r}$, we have
$$\varphi(\hat y_{\lambda_r})\int_T\gamma(t)d\mu\leq
\int_T\gamma(t)\varphi(u(t))d\mu\ ,$$
and hence
$$\varphi(\hat y_{\lambda_r})\int_T\gamma(t)d\mu\leq
\inf_{u\in V_{\gamma,r}}\int_T\gamma(t)\varphi(u(t))d\mu\ .\eqno{(4)}$$
In view of $(3)$, to get $(2)$, we have to show that
$$\varphi(\hat y_{\lambda_r})\int_T\gamma(t)d\mu=
\inf_{u\in V_{\gamma,r}}\int_T\gamma(t)\varphi(u(t))d\mu\ .\eqno{(5)}$$
Arguing by contradiction, assume that $(5)$ does not hold. So, in
view of $(4)$, we would have
$$\varphi(\hat y_{\lambda_r})\int_T\gamma(t)d\mu<
\inf_{u\in V_{\gamma,r}}\int_T\gamma(t)\varphi(u(t))d\mu\ .\eqno{(6)}$$
 From $(6)$, in turn, as
 $(T,{\cal F},\mu)$ is $\sigma$-finite, it would follow the existence of
$\tilde T\in {\cal F}$, with $\mu(\tilde T)<+\infty$, such that
$$\varphi(\hat y_{\lambda_r})\int_{\tilde T}\gamma(t)d\mu<
\inf_{u\in V_{\gamma,r}}\int_T\gamma(t)\varphi(u(t))d\mu\ .\eqno{(7)}$$
Now, consider the function $\hat u:T\to Y$ defined by
$$\hat u(t)=\cases {\hat y_{\lambda_r} & if $x\in \tilde T$\cr & \cr
0 & if $x\in T\setminus\tilde T$\ .\cr}$$
Clearly, $\hat u\in L^p(T,Y)$. We also have
$$\int_T\gamma(t)\psi(\hat u(t))d\mu=
\int_{\tilde T}\gamma(t)\psi(\hat u(t))d\mu\leq r\int_T\gamma(t)d\mu$$
and so $\hat u\in V_{\gamma,r}$. But
$$\int_T\gamma(t)\varphi(\hat u(t))d\mu=\varphi(\hat y_{\lambda_r})
\int_{\tilde T}\gamma(t)d\mu$$
and this contradicts $(7)$. The proof is complete.\hfill $\bigtriangleup$\par
\medskip
REMARK 1. -  In general, the conclusion of Theorem 1 is no longer true if,
for some $\lambda\in ]a,b[$,
the function $\varphi+\lambda\psi$
has more than one global minimum. A simple example (with $a=0$ and $b=+\infty$)
is provided by taking
$Y={\bf R}$,
$$\varphi(y)=\cases {y^2 & if $y\leq 1$\cr & \cr
2-y & if $y>1$\cr}$$
and
$$\psi(y)=y^2\ .$$
So, $\varphi$ is unbounded below and $\varphi+\lambda\psi$ is coercive for all $\lambda>0$.
Clearly, we have $\alpha=0$ and $\beta=+\infty$.
However, for $r=1$, $(2)$ is not satisfied, since $0\in V_{\gamma,r}$,
$\int_T\gamma(t)\varphi(0)d\mu=0$, while $\inf_{\psi^{-1}(1)}\varphi=1$.\par
\medskip
REMARK 2. - At present, we do not know if the conclusion of Theorem 1 holds 
without the coercivity assumption on $\varphi+\lambda\psi$.
\par
\medskip
We now consider a series of consequences of Theorem 1.\par
\smallskip
First, we want to state explicitly the form that Theorem 1 assumes when
$T={\bf N}$, ${\cal F}$ is the power set of ${\bf N}$ and
$$\mu(A)=\hbox {\rm card}(A)$$
for all $A\subseteq {\bf N}$.\par
\smallskip
  Denote by $l_p(Y)$ the
space of all sequences $\{u_n\}$ in $Y$ such that
$$\sum_{n=1}^{\infty}\|u_n\|^p<+\infty\ .$$
\medskip
THEOREM 2. - {\it Let $\varphi, \psi$ satisfy the assumptions of Theorem 1.
\par
Then, for each sequence $\{a_n\}\in l_1({\bf R})\setminus
\{0\}$, with $\inf_{n\in {\bf N}}
a_n\geq 0$,
and for each $r\in ]\alpha,\beta[$, if we put
$$V_{\{a_n\},r}=
\left \{\{u_n\}\in l_p(Y) : \sum_{n=1}^{\infty}a_n\psi(u_n)\leq
r\sum_{n=1}^{\infty}a_n\right\}\ ,$$
we have
$$\inf_{\{u_n\}\in V_{\{a_n\},r}}
\sum_{n=1}^{\infty}a_n\varphi(u_n)=
\inf_{\psi^{-1}(r)}\varphi\sum_{n=1}^{\infty}a_n\ .$$}
\medskip
The next two results deals with consequences of Theorem 1 in the case where
$\varphi\in Y^{*}\setminus \{0\}$.\par
\medskip
THEOREM 3. - {\it Let $y\to \|y\|^q$ be strictly convex for some $q>1$
 and let $\varphi$ be non-zero, continuous and linear.
Moreover, let $\eta:[0,+\infty[\to {\bf R}$ be an increasing strictly
convex function.\par
Then, for each
$\gamma\in L^{\infty}(T)\cap L^1(T)\setminus \{0\}$, with $\gamma\geq 0$,
and for each $r>\eta(0)$ and $p\geq 1$, if we put
$$V_{\gamma,r}=
\left \{u\in L^p(T,Y) : \int_T\gamma(t)\eta(\|u(t)\|^q)d\mu\leq
r\int_T\gamma(t)d\mu\right
\}\ ,$$
we have
$$\inf_{u\in V_{\gamma,r}}
\int_T\gamma(t)\varphi(u(t))d\mu=
-\|\varphi\|_{Y^*}(\eta^{-1}(r))^{1\over q}\int_T\gamma(t)d\mu\ .$$}
\smallskip
PROOF.
 By the assumptions made on $\eta$, the functional
 $y\to \eta(\|y\|^q)$ is strictly convex and, for some $m, c>0$, one
 has
 $$\eta(t)\geq mt-c$$
 for all $t\geq 0$. As a consequence, for each $\lambda>0$, the functional
$y\to \varphi(y)+\lambda\eta(\|y\|^q)$ is coercive and has a unique
global minimum in $X$. At this point, the conclusion follows directly
from Theorem 1, applied taking $a=0$, $b=+\infty$,
$\psi(y)=\eta(\|y\|^q)$ and observing that
$(1)$ holds for each $p\geq 1$ and that $\alpha=\eta(0)$ , $\beta=+\infty$.
\hfill $\bigtriangleup$\par
\medskip
THEOREM 4. - {\it Let $\varphi$ be non-zero, continuous and linear and
let $\psi$ be $C^1$ with
$$\lim_{\|y\|\to +\infty}{{\psi(y)}\over {\|y\|}}=+\infty\ .\eqno{(8)}$$
Finally, assume that, for each $\mu<0$, the equation
$$\psi'(y)=\mu\varphi \eqno{(9)}$$ 
has a unique solution in $Y$ or even at most two when $\hbox {\rm dim}(Y)<
\infty$\ .\par
Then, for each $p\geq 1$, the conclusion of Theorem 1 holds with
any $r>\inf_Y\psi$\ .}\par
\smallskip
PROOF.  In view of $(8)$, the functional $\varphi+\lambda
\psi$ is coercive for each $\lambda>0$. Let $\hat x$ be a global minimum
of this functional. Then, $\hat x$ satisfies $(9)$ with $\mu=-\lambda^{-1}$.
So, when dim$(Y)=\infty$, the uniqueness of $\hat x$ follows from an
assumption directly. Now, assume that dim$(Y)<\infty$. In this case,
$\varphi+\lambda\psi$ satisfies the Palais-Smale condition. As a consequence,
if $\varphi+\lambda\psi$ was admitting two global minima, then,
thanks to Corollary 1 of [3], $(9)$ would have at least three solutions
for $\mu=-\lambda^{-1}$, against an assumption. Now, we can apply
Theorem 1, with $p\geq 1$, $a=0$, $b=+\infty$, observing that
$\alpha=\inf_Y\psi$ and $\beta=+\infty$. \hfill
$\bigtriangleup$\par
\medskip
Here is a consequence of Theorem 1 in the case when $Y$ is a Hilbert
space and $\varphi$ has a Lipschitzian derivative:\par
\medskip
THEOREM 5. - {\it Let $Y$ be a Hilbert space,
let $\varphi$ be $C^1$ and let $\varphi'$ be Lipschitzian,
with Lipschitz constant $L>0$. Assume that $\varphi'(0)\neq 0$.
Set
$$S=\{y\in Y: \varphi'(y)+L y=0\}$$
and
$$\rho=\inf_{y\in S}\|y\|^2\ .$$
Then, for each
$\gamma\in L^{\infty}(T)\cap L^1(T)\setminus \{0\}$,
with $\gamma\geq 0$,
and for each $r\in ]0,\rho[$, $p\geq 2$, if we put
$$V_{\gamma,r}=
\left \{u\in L^p(T,Y) : \int_T\gamma(t)\|u(t)\|^2d\mu\leq
r\int_T\gamma(t)d\mu\right
\}\ ,$$
we have
$$\inf_{u\in V_{\gamma,r}}\int_T\gamma(t)
\varphi(u(t))d\mu=\inf_{\|y\|^2=r}\varphi(y)\int_T\gamma(t)d\mu
\ .$$}
\smallskip
PROOF. 
Note that the functional $y\to \varphi(y)+{{\lambda}\over {2}}\|y\|^2$
is convex if $\lambda=L$, while it is strictly convex and coercive
if $\lambda>L$ 
(see, for instance, Proposition 2.2 of [6]). So, this functional
has a unique global minimum if
$\lambda>L$, while the set of its global minima coincides with $S$ if
$\lambda=L$.
At this point, the conclusion is obtained applying Theorem 1
with
$$\psi(y)={{\|y\|^2}\over {2}}$$
for all $y\in Y$ and
$$a=L\ , b=+\infty\ ,$$
taking into account that $(1)$ is satisfied for each $p\geq 2$ since
$\varphi'$ is Lipschitzian and observing that $\alpha=0$ and $\beta={{\rho}\over
{2}}$.
\hfill $\bigtriangleup$
\medskip
REMARK 3. - Note that Theorem 5 is an extension of Theorem 1 of [7].\par
\medskip
In the next result, we will apply Theorem 1 taking as $Y$ the usual
Sobolev space $W^{1,q}_0(\Omega)$ with the usual norm
$$\left ( \int_{\Omega}|\nabla v(x)|^qdx\right )^{1\over q} ,$$
where $\Omega$ is bounded domain in ${\bf R}^n$ ($n\geq 3$)
 with smooth boundary and $q>1$.\par
\smallskip
Moreover, if $u\in L^p(T,W^{1,q}_0(\Omega))$ we will write $u(t,x)$ instead
of $u(t)(x)$. That is, we will identify $u$ with
the function $\omega:T\times \Omega\to {\bf R}$ defined by
$$\omega(t,x)=u(t)(x)$$
for all $(t,x)\in T\times\Omega$.  \par
\medskip
THEOREM 6. - {\it Let $f:{\bf R}\to [0,+\infty[$ be a continuous
function, with $f(0)=0$ and
$\liminf_{\xi\to +\infty}f(\xi)>0$,
 such that $\xi\to {{f(\xi)}\over {\xi^{q-1}}}$
is decreasing in $]0,+\infty[$ and
$$\lim_{|\xi|\to +\infty}{{f(\xi)}\over {|\xi|^{q-1}}}=0 \eqno{(10)}$$
for some $q>1$.\par
Then, for each
$\gamma\in L^{\infty}(T)\cap L^1(T)\setminus \{0\}$, with $\gamma\geq 0$,
and each $r>0$, $p\geq q$, if we put
$$V_{\gamma,r}=\left \{u\in L^p(T,W^{1,q}_0(\Omega)) : \int_T\gamma(t)\left (
\int_{\Omega}|\nabla u(t,x)|^qdx\right ) d\mu\leq r\int_T\gamma(t)d\mu\right \}\ ,$$
we have
$$\sup_{u\in V_{\gamma,r}}\int_T\gamma(t)\left ( \int_{\Omega}F(u(t,x))dx\right )d\mu=
\sup_{v\in W^{1,q}_0(\Omega), \int_{\Omega}|\nabla v(x)|^qdx=r}
\int_{\Omega}F(v(x))dx\int_T\gamma(t)d\mu\ ,$$
where 
$$F(\xi)=\int_0^{\xi}f(s)ds$$
for all $\xi\in {\bf R}$.} \par
\smallskip
PROOF. We are going to apply Theorem 1 taking $Y=W^{1,q}_0(\Omega)$ and
$$\varphi(v)=-\int_{\Omega}F(v(x))dx\ ,$$
$$\psi(v)=\int_{\Omega}|\nabla v(x)|^qdx$$
for all $v\in W^{1,q}_0(\Omega)$. Due to $(10)$,
 by classical results, 
$\varphi$ is sequentially weakly continuous in $W^{1,q}_0(\Omega)$,
$(1)$ is satisfied for any $p\geq q$, and, for each $\lambda>0$,
the functional $\varphi+\lambda\psi$ is $C^1$, coercive and
satisfies the Palais-Smale condition. Moreover, since $f\geq 0$,
its non-zero critical points are strictly positive in $\Omega$ ([1], [8]).
 Moreover, since
the function $\xi\to {{f(\xi)}\over {\xi^{q-1}}}$
is decreasing in $]0,+\infty[$, Proposition 4.2 of [2] ensures that,
for each $\lambda>0$, there exists at most one strictly
positive critical point of $\varphi+\lambda\psi$. As a consequence,
we infer that, for each $\lambda>0$, the functional $\varphi+\lambda\psi$
has a unique global minimum in $W^{1,q}_0(\Omega)$, since otherwise,
in view of Corollary 1 of [3], it would have at least three critical
points. Hence, we are allowed to apply Theorem 1 with $a=0$ and
$b=+\infty$. Clearly, we have $\alpha=0$ and $\beta=+\infty$ (since
$\lim_{\xi\to +\infty}F(\xi)=+\infty$ and hence $\varphi$ is
unbounded below). The proof is complete.\hfill $\bigtriangleup$\par
\medskip
The next application of Theorem 1 concerns
a Jensen-like inequality in $L^p$-spaces.\par
\medskip
THEOREM 7. - {\it Let $f:{\bf R}\to {\bf R}$
 be a continuous function, differentiable in $]0,+\infty[$,
with $\sup_{]-\infty,0]}f\leq 0$. Assume that, for some
$\delta\geq 0$, the
function $y\to \delta|y|^p-f(y)$ has no global minima in
${\bf R}$, 
$$\limsup_{y\to +\infty}{{f(y)}\over {y^p}}=\delta \eqno{(11)}$$
and the function
$$y\to {{f'(y)}\over {y^{p-1}}}$$
 is injective in $]0,+\infty[$.\par
Then, for each $\gamma\in L^{\infty}(T)\cap L^1(T)\setminus \{0\}$,
with $\gamma\geq 0$, 
one has
$$\int_T\gamma(t)f(u(t))d\mu\leq
f\left (\left ( {{\int_T\gamma(t)|u(t)|^pd\mu}\over {\int_T\gamma(t)d\mu}}\right )
^{1\over p}\right )\int_T\gamma(t)d\mu\ ,$$
for all $u\in L^p(T)$\ .}
\smallskip
PROOF. We are going to apply Theorem 1 with $Y={\bf R}$, $\varphi(y)=
-f(y)$, $\psi(y)=|y|^p$ and $a=\delta$, $b=+\infty$. 
 Fix
$\lambda>\delta$.
 From $(11)$, we clearly infer that $\varphi+\lambda\psi$ is coercive. We
now show that this function has a unique global minimum. Arguing by
contradiction, assume that $y_1, y_2\in {\bf R}$ are two distinct
global minima of $\varphi+\lambda\psi$. We can suppose that $y_1<y_2$.
Since $\varphi(y)+\lambda\psi(y)>0$ for all $y<0$ and $\varphi(0)+\lambda\psi(0)=0$,
it would follow that $y_1\geq 0$. By the Rolle theorem, there would be
$y_3\in ]y_1,y_2[$ such that 
$$p\lambda y_3^{p-1}=f'(y_3)\ .$$
As a consequence, we would have
$${{f'(y_2)}\over {y_2^{p-1}}}={{f'(y_3)}\over {y_3^{p-1}}}\ ,$$
contrary to the assumption that the function $y\to {{f'(y)}\over {y^{p-1}}}$ is
injective in $]0,+\infty[$.
So, we are allowed to apply Theorem 1, observing that
$\alpha=0$ and $\beta=+\infty$. Let $u\in L^p(T)\setminus \{0\}$. Put
$$r={{\int_T\gamma(t)|u(t)|^pd\mu}\over {\int_T\gamma(t)d\mu}}\ .$$
Clearly, we have
$$\inf_{\psi^{-1}(r)}\varphi=-f\left (\left ( {{\int_T\gamma(t)|u(t)|^pd\mu}\over {\int_T\gamma(t)d\mu}}\right )
^{1\over p}\right )$$
and hence, since $u\in V_{\gamma,r}$, it follows
$$\int_T\gamma(t)f(u(t))d\mu\leq
f\left (\left ( {{\int_T\gamma(t)|u(t)|^pd\mu}\over {\int_T\gamma(t)d\mu}}\right )
^{1\over p}\right )\int_T\gamma(t)d\mu\ ,$$
as claimed.\hfill$\bigtriangleup$
\medskip
REMARK 4. -
The class of functions $f$ satisfying the assumptions of Theorem 7 is quite
broad. For instance, a typical function in that class is
$$f(y)=a_0\log(1+(y^+)^p)+\sum_{i=1}^k a_i (y^+)^{q_i}$$
where $y^+=\max\{y,0\}$,
$a_i$ ($i=0,...,k$) are $k+1$ non-negative numbers, with
$\sum_{i=0}^k a_i>0$, and $q_i$  ($i=1,...,k$)
are $k$ positive numbers less than $p$.
\medskip
As a consequence of this remark, we get, for instance, the following\par
\medskip
COROLLARY 1. - {\it For each $\gamma\in L^{\infty}(T)\cap L^1(T)\setminus
\{0\}$, with $\gamma\geq 0$, one has
$$\int_T\gamma(t)\log(1+(u(t))^p)d\mu\leq
\log\left ( 1+{{\int_T\gamma(t)(u(t))^pd\mu}\over {\int_T\gamma(t)
d\mu}}
\right )\int_T\gamma(t)d\mu \eqno{(12)}$$
for all $u\in L^p(T)$ with $u\geq 0$.}
\medskip
Assume that $\mu(T)=\gamma=1$. It is worth noticing that, in this case,
 $(12)$
 can be obtained by the classical Jensen
inequality only when $p=1$. In fact, when $p>1$, the function $t\to
\log(1+t^p)$ is neither concave nor convex in $[0,+\infty[$. While, when
$p<1$, the use of the Jensen inequality would provide
$$\int_T\log(1+(u(t))^p)d\mu\leq
\log\left (1+\left (\int_Tu(t)d\mu\right ) ^p\right ) \ .$$ 
Note that this latter inequality is weaker than $(12)$ since
$$\int_T(u(t))^pd\mu\leq
\left (\int_Tu(t)d\mu\right ) ^p\ .$$
\medskip
The final result is an application of Theorem 1 to quasi-linear
equations.\par
\smallskip
So, in the sequel, $\Omega\subseteq {\bf R}^n$ is a bounded domain
with smooth boundary and $p>1$.\par
\smallskip
If $n\geq p$, we denote by ${\cal A}$ the class of all
continuous functions $f:{\bf R}\to {\bf R}$ such that
$$\sup_{y\in {\bf R}}{{|f(y)|}\over
{1+|y|^{s}}}<+\infty\ ,$$
where  $0<s< {{pn-n+p}\over {n-p}}$ if $p<n$ and $0<s<+\infty$ if
$p=n$. While, when $n<p$, ${\cal A}$ stands for the class
of all continuous functions $f:{\bf R}\to {\bf R}$. 
Given $f\in {\cal A}$, consider the following Dirichlet problem
$$\cases {-\hbox {\rm div}(|\nabla u|^{p-2}\nabla u)=
f(u)
 & in
$\Omega$\cr & \cr u=0 & on
$\partial \Omega$\ .\cr}\eqno{(P_{f})} $$
 Let us recall
that a weak solution
of $(P_{f})$ is any $u\in W^{1,p}_0(\Omega)$ such that
 $$\int_{\Omega}|\nabla u(x)|^{p-2}\nabla u(x)\nabla v(x)dx
-\int_{\Omega}
f(u(x))v(x)dx=0$$
for all $v\in W^{1,p}_0(\Omega)$.\par
\smallskip
Moreover, $\lambda_{1,p}$ denotes the principal
eigenvalue of the problem
$$\cases {-\hbox {\rm div}(|\nabla u|^{p-2}\nabla u)=
\lambda  |u|^{p-2}u
 & in
$\Omega$\cr & \cr u=0 & on
$\partial \Omega$\ .\cr} $$
We have
$$\lambda_{1,p}=\inf_{u\in W^{1,p}_0(\Omega)\setminus \{0\}}{{\int_{\Omega}|\nabla u(x)|^pdx}\over
{\int_{\Omega}|u(x)|^pdx}}\ .$$
\medskip
Also, let us recall the following consequence of the variational principle established in [4]:\par
\medskip
THEOREM B. - {\it Let $X$ be a reflexive real Banach space and let
$\Phi, \Psi:X\to {\bf R}$ be two sequentially weakly lower semicontinuous
functionals, with $\Phi(0)=\Psi(0)=0$, and with $\Psi$ also coercive and continuous.\par
Then, for each $\sigma>\inf_X\Psi$ and each $\lambda$ satisfying
$$\lambda>-{{\inf_{\Psi^{-1}(]-\infty,\sigma])}\Phi}\over {\sigma}}$$
the functional $\lambda\Psi+\Phi$  has a local
minimum belonging to $\Psi^{-1}(]-\infty,\sigma[)$\ .}\par
\medskip
The final result is as follows:\par
\medskip
THEOREM 8. - {\it Let $f\in {\cal A}$, with $f\geq 0$, 
 and let $F(y)=\int_0^yF(t)dt$ for all $y\in {\bf R}$. 
Assume that:\par
\noindent
$(a_1)$\hskip 5pt $\lim_{y\to 0^+}{{F(y)}\over {y^p}}=+\infty\ ;$\par
\noindent
$(a_2)$\hskip 5pt $\delta:=\limsup_{y\to +\infty}{{F(y)}\over {y^p}}<+\infty\ ;$\par
\noindent
$(a_3)$\hskip 5pt the function $y\to \delta y^p-F(y)$ has no global minima in $[0,+\infty[$\ ;\par
\noindent
$(a_4)$\hskip 5pt 
for each $\lambda>p\delta$, the equation $\lambda y^{p-1}=f(y)$ has
at most two solutions in $]0,+\infty[$\ .\par
\noindent
Under such hypotheses, for each $\rho>0$ and each $\nu\in ]0,1]$ satisfying
$$\nu<{{\lambda_{1,p}\rho^p}\over {p F(\rho)}}\ ,\eqno{(13)}$$
  the problem
$$\cases {-\hbox {\rm div}(|\nabla u|^{p-2}\nabla u)=
\nu f(u)
 & in
$\Omega$\cr & \cr u=0 & on
$\partial \Omega$\cr} $$
has a positive weak solution satisfying
$$\int_{\Omega}|\nabla u(x)|^pdx<\rho^p\lambda_{1,p}\hbox {\rm meas}(\Omega)\ .$$}
\smallskip
PROOF. Fix $\rho$ and $\nu$ as above.  
Since $f\geq 0$, by classical results ([1], [8]), the positive weak solutions of
the problem are exactly the non-zero critical points in $W^{1,p}_0(\Omega)$ of the
energy functional 
$$u\to {{1}\over {p}}\int_{\Omega}|\nabla u(x)|^pdx-\nu\int_{\Omega}F(u(x))dx\ .$$
 We are going to apply Theorem 1 taking $Y={\bf R}$, $\varphi(y)=-\nu F(y)$, $\psi(y)=|y|^p$,
$a=\delta$ and $b=+\infty$. Note that $\varphi$ is
non-negative in $]-\infty,0]$. So, $(1)$ is satisfied in view of $(a_2)$.
Fix $\lambda>\delta$.
From $(a_2)$ again, it follows that $\varphi+\lambda\psi$ is coercive . Arguing by contradiction,
assume that $\varphi+\lambda\psi$ has two global minima, say $y_1, y_2$, with
$y_1<y_2$. Differently from Theorem 7, this time we are assuming $(a_1)$ from which it
follows that 
$$\inf_{[0,+\infty[}(\varphi+\lambda\psi)<0\ .$$
This fact implies that $y_1>0$. As a consequence,
 the equation
$$p\lambda y^{p-1}=\nu f(y)$$
would admit the solutions $y_1$, $y_2$ and a third one in $]y_1,y_2[$ given by
the Rolle theorem. But, this contradicts $(a_4)$ .
 Hence, the function
$\varphi+\lambda\psi$ has a unique global minimum. Further, note that $\alpha=0$ 
and, in view of $(a_3)$, $\beta=+\infty$. Then, if we put
$$V_{\rho}=\left\{u\in L^p(\Omega) : \int_{\Omega}|u(x)|^pdx\leq \rho^p\hbox {\rm meas}(\Omega)\right\}\ ,$$
Theorem 1 ensures that
$$\sup_{u\in V_{\rho}}\int_{\Omega}F(u(x))dx=F(\rho)\hbox {\rm meas}(\Omega)\ .\eqno{(14)}$$
On the other hand, setting
$$B_{\rho}=\left\{u\in W^{1,p}_0(\Omega) : \int_{\Omega}|\nabla u(x)|^pdx\leq 
\rho^p\lambda_{1,p}\hbox {\rm meas}(\Omega)\right\}\ ,$$
we have
$$B_{\rho}\subseteq V_{\rho}\ .$$
Consequently
$$\sup_{u\in B_{\rho}}\int_{\Omega}F(u(x))dx\leq \sup_{u\in V_{\rho}}\int_{\Omega}F(u(x))dx\ .\eqno{(15)}$$
Now, if we put
$$\sigma=\rho^p\lambda_{1,p}\hbox {\rm meas}(\Omega)\ ,$$
in view of $(13)$ , $(14)$ and $(15)$, we have
$$\sup_{u\in W^{1,p}_0(\Omega), \int_{\Omega}|\nabla u(x)|^pdx\leq \sigma}
\int_{\Omega}\nu F(u(x))dx<{{\sigma}\over {p}}\ .$$
At this point, we can apply Theorem B taking $X=W^{1,p}_0(\Omega)$,
$\Psi(u)=\int_{\Omega}|\nabla u(x)|^pdx$ and $\Phi(u)=-\nu\int_{\Omega}F(u(x))dx$.
Hence, the energy functional has a local minimum  $u$ (which is therefore a solution of the problem)  such that 
$$\int_{\Omega}|\nabla u(x)|^pdx<\rho^p\lambda_{1,p}\hbox {\rm meas}(\Omega)\ .$$
To show that $u\neq 0$, we finally remark that $0$ is not a local minimum of the
energy functional. Indeed, by a classical result, there is a
bounded and positive function $v\in W^{1,p}_0(\Omega)$ such that
$$\int_{\Omega}|\nabla v(x)|^pdx=\lambda_{1,p}\int_{\Omega}|v(x)|^pdx\ .$$
By $(a_1)$, there is $\theta>0$ such that
$$F(y)>{{\lambda_{1,p}}\over {\nu p}}y^p$$
for all $y\in ]0,\theta[$. Hence, for each $\eta\in \left ] 0, {{\theta}\over {\sup_{\Omega}v}}\right[$,
we have
$$\nu\int_{\Omega}F(\eta v(x))dx>{{\lambda_{1,p}}\over {p}}\int_{\Omega}|\eta v(x)|^pdx=
{{1}\over {p}}\int_{\Omega}|\nabla \eta v(x)|^pdx\ .$$ 
This shows that the energy functional takes negative values in each ball of $W^{1,p}_0(\Omega)$
centered at $0$ and so $0$ is not a local minimum for it. The proof is complete. \hfill
$\bigtriangleup$\par
\medskip
Note the following corollary of Theorem 8 (for the uniqueness, consider again Proposition 4.2 of [2]):\par
\medskip
COROLLARY 2. - {\it For each $\nu\in ]0,1]$, the unique positive weak solution of the problem
$$\cases {-\hbox {\rm div}(|\nabla u|\nabla u)=
\nu u
 & in
$\Omega$\cr & \cr u=0 & on
$\partial \Omega$\cr} $$
satisfies the inequality
$$\int_{\Omega}|\nabla u(x)|^3dx\leq {{27\hbox {\rm meas}(\Omega)}\over {8\lambda_{1,3}^2}}\nu^3\ .$$}
\bigskip
\bigskip
\vfill\eject
\centerline {\bf References}\par
\bigskip
\bigskip
\noindent
[1]\hskip 5pt P. DR\'ABEK, {\it On a maximum principle for weak solutions of some
quasi-linear elliptic equations}, Appl. Math. Letters, {\bf 22}
(2009), 1567-1570.\par
\smallskip
\noindent
[2]\hskip 5pt P. DR\'ABEK and J. HERN\'ANDEZ, {\it Existence and
uniqueness of positive solutions for some quasilinear elliptic problems},
Nonlinear Anal., {\bf 44} (2001), 189-204.\par
\smallskip
\noindent
[3]\hskip 5pt  P. PUCCI and J. SERRIN, {\it A mountain pass theorem},
J. Differential Equations, {\bf 60} (1985), 142-149.\par
\smallskip
\noindent
[4]\hskip 5pt B. RICCERI, {\it A general variational principle and
some of its applications}, J. Comput. Appl. Math., {\bf 113}
(2000), 401-410.\par
\smallskip
\noindent
[5]\hskip 5pt B. RICCERI, {\it Well-posedness of constrained minimization
problems via saddle-points}, J. Global Optim., {\bf 40} (2008),
389-397.\par
\smallskip
\noindent
[6]\hskip 5pt B. RICCERI, {\it Fixed points of nonexpansive potential operators in Hilbert
spaces}, Fixed Point Theory Appl., {\bf 2012}, 2012: 123. \par
\smallskip
\noindent
[7]\hskip 5pt B. RICCERI, {\it Another fixed point theorem for nonexpansive
potential operators}, Studia Math., {\bf 211} (2012), 147-151.\par
\smallskip
\noindent
[8]\hskip 5pt J. L. V\'AZQUEZ, {\it A strong maximum principle for some
quasilinear elliptic equations}, Appl. Math. Optim., {\bf 12} (1984),
191-202. \par
\bigskip
\bigskip
\bigskip
\bigskip
Department of Mathematics\par
University of Catania\par
Viale A. Doria 6\par
95125 Catania\par
Italy\par
{\it e-mail address}: ricceri@dmi.unict.it

\bye